\RequirePackage{fix-cm}
\documentclass[11pt]{article}
\usepackage{graphicx,amsmath,amsfonts,amscd,amssymb,bm,cite,epsfig,epsf,url,color,algorithm,algorithmic,enumerate}
\usepackage{fullpage}
\usepackage[small,bf]{caption}
\usepackage{subfigure}
\usepackage{lscape}
\usepackage[space]{grffile}
\usepackage{epsfig,psfrag,amstext,subfigure,subeqnarray,nccfoots,color,multirow,bm}
\usepackage{appendix}
\newtheorem{theorem}{Theorem}[section]
\newtheorem{lemma}[theorem]{Lemma}

\newtheorem{proposition}[theorem]{Proposition}
\newtheorem{definition}[theorem]{Definition}

\newtheorem{assumption}[theorem]{Assumption}

\newcommand{\dist}{{\rm dist}}
\newcommand{\smax}{\sigma_{\max}}

\numberwithin{equation}{section}
\def \endprf{\hfill {\vrule height6pt width6pt depth0pt}\medskip}
\newenvironment{proof}{\noindent {\bf Proof} }{\endprf\par}

\title{\bf A global  dual error bound and its application to the analysis of linearly constrained nonconvex optimization\thanks{This research is supported in part by the NSFC grants 61731018 (key project) and 61571384, and by the Peacock project of Shenzhen Municipal Government. }}
\author{Jiawei Zhang\thanks{Shenzhen Research Institute of Big Data, The Chinese University of Hong Kong, Shenzhen, China.
Email: luozq@cuhk.edu.cn}
\ and \
Zhi-Quan Luo$^\dagger$}
\begin{document}
\maketitle
\begin{abstract}
Error bound analysis, which estimates the distance of a point to the solution set of an optimization problem using the optimality residual, is a powerful tool for the  analysis of first-order optimization algorithms.
In this paper, we use global error bound analysis to study the iteration complexity of a first-order algorithm for a linearly constrained nonconvex minimization problem.
we develop a global dual error bound analysis for a regularized version of this nonconvex problem by using  a novel ``decomposition'' technique.
Equipped with this global dual error bound, we prove that a suitably designed primal-dual first order method can generate an $\epsilon$-stationary solution of the linearly constrained nonconvex minimization problem within $\mathcal{O}(1/\epsilon^2)$ iterations, which is the best known iteration complexity for this class of nonconvex problems.

\end{abstract}
\section{Introduction}
Consider the following linearly constrained optimization problem:
\begin{equation}\label{P1}
\begin{array}{ll}
\mbox{minimize}& f(x)\\ [5pt]
\mbox{subject to} & Ax=b,\ x\in P,
\end{array}
\end{equation}
where $P=\{x\mid Gx\preceq h\}$ is a polyhedral set in $\mathbb{R}^n$, $A\in \mathbb{R}^{m\times n}$ and $G\in \mathbb{R}^{l\times n}$ are  two given matrices. The objective function $f$ is assumed to be smooth but possibly nonconvex, whose gradient is Lipschitz-continuous (with a  Lipschitz constant $L_f$).
Also, we assume that $f$ is bounded from below over the feasible set $\{x\mid Ax=b, x\in P\}$.
Problem \eqref{P1} appears in many practical applications, such as resource allocation \cite{yan2020collaborative}, rate maximization in wireless communication \cite{scutari2008optimal}, clustering \cite{bauckhage2015k,ding2008equivalence}, non-negative matrix factorization  \cite{hajinezhad2016nonnegative} and distributed optimization \cite{Hong17,zeng2018nonconvex}.

A popular approach to solve \eqref{P1} is using a primal-dual first order method which alternately updates primal and dual iterates via inexact gradient steps  \cite{Hong17,tebulle}. Such methods are well suited for large scale optimization problems involving big data. However, the convergence analysis of this type of primal-dual first order methods is a well known difficult problem especially if $f$ is nonconvex. One particularly effective technique to tackle this problem is to use error bound analysis \cite{Pang-survey, eb-survey, Pang-complementary}) which uses the optimality residuals of \eqref{P1} to estimate the distance of a point to the solution set of \eqref{P1}. This approach has been successful in the convergence rate analysis of first-order methods for convex $f$ \cite{Luo-dual,bcd,hong2017linear} as well as in a nonconvex setting \cite{zhang2018proximal}. One weakness of these works is that the convergence rate is ``local" (meaning linear convergence is only guaranteed when the iterates are close to the solution set), and dependent on the error bound constant which is either unknown or difficult to estimate. The ``localness" of the convergence result is due to the error bound being only shown to hold locally around the solution set of the problem.

\subsection{Related works}
There are many studies of error bound analysis in the optimization literature \cite{eb-survey,Pang-survey}. For instance, in \cite{Luo-dual}, authors prove the ``local'' linear convergence of the dual ascent algorithm  for a family of convex problems with polyhedral constraints in the absence of strong convexity. In reference \cite{bcd}, authors use ``local'' error bounds to show the linear convergence of block coordinate descent for a family of convex problem, while in \cite{hong2017linear}, authors show the linear convergence of ADMM algorithm using an error bound approach, again without strong convexity.

Recently, \cite{zhang2018proximal} used ``local'' dual error bound to deal with nonconvex problems and prove the convergence of a smoothed proximal augmented Lagrangian method. Papers such as \cite{So19} also use the error bound  to analyze second-order optimization methods.

For nonconvex problem, there are some recent papers focusing on iteration complexity of first-order  algorithms. In \cite{kong2018complexity}, the authors propose  a quadratic penalty accelerated inexact proximal point method with $O(1/\epsilon^3)$ iteration complexity for finding an $\epsilon$-stationary solution of a linearly constrained nonconvex composite problem. Recently, the authors of \cite{lin2019inexact} propose an inexact proximal-point penalty method for constrained optimization problems, with nonconvex composite objective function and convex constraints, achieving an iteration complexity $\mathcal{O}(1/\epsilon^{2.5})$, which is by far the best result for nonconvex problems with convex constraints.
In contrast to the double-loop method \cite{kong2018complexity,lin2019inexact}, a single-loop perturbed proximal primal-dual algorithm is proposed  in \cite{hajinezhad2019perturbed}, and the authors  analyze its asymptotic convergence, though no iteration complexity analysis is given.

\subsection{Our contributions}

Our contributions are as follows:

\begin{enumerate}
  \item We develop a novel technique to establish a global dual error bound for a regularized version of problem \eqref{P1}. This technique allows us to estimate the distance to the solution set from any point rather than just points near the solution set. Compared to \cite{zhang2018proximal}, we remove the slater condition, strict complementarity assumption and the compactness assumption of the feasible set.
  \item Using the global dual error bound analysis, we can explicitly compute the parameters of the smoothed prox-ALM algorithm proposed in\cite{zhang2018proximal}. We also show that the iteration complexity of the algorithm is $B/\epsilon^2$, where $B>0$ is a global constant. Compared to the iteration complexity in \cite{lin2019inexact}, our complexity bound is lower by a factor of $1/\sqrt{\epsilon}$ for linearly constrained problems. Though the results of \cite{lin2019inexact} are for a more general setting,  problem \eqref{P1} is an important special case and our iteration complexity achieves the optimal order $O(1/\epsilon^2)$ for this type of nonconvex optimization problems.
\end{enumerate}

\noindent{\bf Remark.} Note that most double-loop algorithms for constrained nonconvex problems are based on the proximal framework, requiring at least $\mathcal{O}(1/\epsilon^2)$ outer iterations and a number of inner iterations that also grows with $1/\epsilon$. Hence the iteration complexity of the double-loop algorithms for nonconvex problems usually need more than $\mathcal{O}(1/\epsilon^2)$ iterations to generate an $\epsilon$-stationary solution.
To achieve the optimal complexity $\mathcal{O}(1/\epsilon^2)$, we need to consider single-loop algorithms with constant step sizes and establish their convergence using global error bounds.

\section{ Preliminaries }
In this section, we give the notations and main definitions. Then we introduce the algorithm discussed in this paper and state our main result.
\subsection{Notations}
First, we list some notations used in this paper as follows.
\begin{enumerate}
\item $[\cdot]_+$ means the projection to the set $P$.
%\item $\sigma$ is the largest singular value of $A$.
\item For a matrix $\mathbf{M}$, $\sigma_{\max}(\mathbf{M})$ and $\sigma_{\min}(\mathbf{M})$
 are the largest and smallest singular values of $\mathbf{M}$ respectively.

 \item $\mathrm{dist}(v, S)$ means the Euclidian distance from a point $v$ to a set $S$.
\item $[l]=\{1, 2, \cdots, l\}$.
\item For a vector $v$,  $v_i$ means the $i$-th component of $v$. For a set $\mathcal{S}$, $v_{\mathcal{S}}\in \mathbb{R}^{|\mathcal{S}|}$ is the vector containing all components  $v_i$'s with $i\in\mathcal{S}$.
\end{enumerate}
\subsection{Our main assumptions}
In this paper, we make the following assumptions:
\begin{assumption}
\label{ass}

\begin{enumerate}
\item [{\rm (a)}] $f$ is a smooth function and $\nabla f(x)$ is $L_f$-Lipschitz continuous, i.e., for any $x, x'\in P$, we have
$$\|\nabla f(x)-\nabla f(x')\|<L_f\|x-x'\|.$$
\item [{\rm (b)}] $f$ is bounded from below in the feasible set $\{x\in P\mid Ax=b\}$, i.e.,
$$f(x)>\underline{f}>-\infty,\quad x\in\{x\in P\mid Ax=b\}$$
for some constant $\underline{f}$.
\end{enumerate}

\end{assumption}
Let $\gamma_f=-L_f$.
Then in view of Assumption~\ref{ass}(a),  $f$ is $\gamma_f$-weakly convex, i.e.,  we have
$$\langle \nabla f(x)-\nabla f(x'), x-x'\rangle\ge \gamma_f\|x-x'\|^2.$$
In other words, for arbitrary fixed $v\in \mathbb{R}^n$, the function $f(x)+\frac{p}{2}\|x-v\|^2$ is a strongly convex function with modulus $(p+\gamma_f)=(p-L_f)$ if $p>L_f$.

\subsection{Stationary solution set of \eqref{P1}}
By the linearity of the constraints, an optimal solution $x^*$ of Problem \eqref{P1} must satisfy the following Karush-Kuhn-Tucker (KKT) conditions:
\begin{eqnarray}\label{KKTforP}
\nabla f(x^*)+A^Ty^*+G^T\mu^*&=&0\\
Ax^*-b&=&0,\\
Gx^*&\preceq& h,\\
\mu^*&\succeq& 0,\\
\mu_i^*(Gx^*-h)_i&=&0,\quad i\in[l],
\end{eqnarray}
where $y^*, \mu^*$ are the Lagrangian multipliers corresponding to the equality and inequality constraints of \eqref{P1}.
 For nonconvex problems, it is NP-hard to get a global minimum of \eqref{P1}.
However, with good initialization, a KKT solution is usually good enough for  practical applications.
In this paper, we focus on finding such a KKT solution of Problem \eqref{P1} using a first-order algorithm.
\begin{definition}
Let  $X^*$ be the stationary solution set of Problem \eqref{P1}. Specifically, $x^*\in X^*$ if there exist some $y^*, \mu^*$ such that the conditions \eqref{KKTforP} hold with $y^*, \mu^*$.
Also let $U^*$ be the primal-dual solution set of \eqref{P1}, i.e., $U^*$ is the set of all pairs $(x^*, y^*)$ satisfying \eqref{KKTforP} with some $\mu^*$.
\end{definition}

%[[[SHOULD DEFINE $\epsilon$-stationary solution and compare with different definitions]]]

Next, we define the $\epsilon$-stationary solution of \eqref{P1} as in \cite{kong2018complexity}. Let $\iota(x)$ be the indicator function of the set $P$, i.e. , $\iota(x)=0$ if $x\in P$ and $\iota(x)=\infty$ otherwise.
\begin{definition}
A primal-dual vector $(x, y)$ is said to be an $\epsilon$-stationary solution of \eqref{P1} if $\|Ax-b\|\le \epsilon$ and  there exists a vector  $v\in \nabla f(x)+A^Ty+\partial{\iota(x)}$ with $\|v\|\le \epsilon$. Here $\partial{\iota(x)}$ is the sub-differential set of $\iota(\cdot)$ at $x$.
\end{definition}

\subsection{Primal-dual first order algorithms for \eqref{P1}}
The augmented Lagrangian function $L_{\rho}(x; y)$ of \eqref{P1} is given by:
$$L_{\rho}(x; y)=f(x)+y^T(Ax-b)+\frac{\rho}{2}\|Ax-b\|^2.$$
The augmented Lagrangian method (ALM) for solving \eqref{P1} is given as follows (see \cite{bertsekas}):
\begin{algorithm}[ht]
	\caption{ALM}
	\label{Alg:ALM}
\begin{algorithmic}[1]
\STATE Let $\rho>0$;
\STATE Initialize $x^0, y^0$;
\FOR{$t=0,1,2,\ldots,$}
\STATE $x^{t+1}=\arg\min_{x\in P}L_{\rho}(x; y^t)$;
\STATE $y^{t+1}=y^t+\rho(Ax^{t+1}-b)$.
\ENDFOR
	%	\Ensure
\end{algorithmic}
\end{algorithm}

The ALM is known to converge when $f$ is convex and satisfies some mild assumptions (\cite{hong2017linear}).
However, the sub-problem
\begin{equation}\label{exactmin}
x^{t+1}=\arg\min_{x\in P}L_{\rho}(x; y^t)
\end{equation}
is usually hard to solve.
Moreover, it is shown by a counter-example in \cite{wyin16} that ALM may not converge if $f$ is nonconvex.
Therefore, a suitable modification for ALM is needed for convergence.
One such modification is to replace the exact minimization step \eqref{exactmin} by a linearized-proximal step which performs a gradient descent to the augmented function, while the dual update is kept unchanged (i.e., still use the constraint residual to update the dual variable). This modified ALM is shown to be convergent if $P=\mathbb{R}^n$ in \cite{Hong17}. The iteration complexity is proved to be $\mathcal{O}(1/\epsilon^2)$.
However, a numerical experiment in \cite{zhang2018proximal} shows that this modified ALM may not converge if $P\ne \mathbb{R}^n$.
Recently, authors of \cite{zhang2018proximal} propose a ``smoothed prox-ALM'', which is proved to converge under some regularity assumptions for the case where $P$ is a bounded box. The iteration complexity is also shown to be $\mathcal{O}(1/\epsilon^2)$.

Let
\begin{equation}\label{eq:K}
K(x, z;y)=L_{\rho}(x;y)+\frac{p}{2}\|x-z\|^2,
\end{equation}
where $p>L_f$ is a positive constant.
Note that $K(x, z; y)$ is strongly convex of $x$ with modulus $\gamma_K=(p+\gamma_f)=(p-L_f)$ and $\nabla_xK(x, z;y)$ is $L_f+\rho\smax^2(A)+p$-Lipschitz-continuous of $x$.
The smoothed prox-ALM algorithm \cite{zhang2018proximal} is given as follows.
\begin{algorithm}[ht]
	\caption{S-prox-ALM}
	\label{Alg:SProxALM}
\begin{algorithmic}[1]
\STATE Let $\rho>0$, $\alpha>0$, $0<\beta\le 1$ and $\frac{1}{L_f+\rho\smax^2(A)+p}>c>0$;
\STATE Initialize $x^0, z^0, y^0$;
\FOR{$t=0,1,2,\ldots,$}
\STATE $y^{t+1}=y^t+\alpha(Ax^t-b)$;
\STATE $x^{t+1}=[x^t-c\nabla_x K(x^t, z^t; y^{t+1})]_+$;
\STATE $z^{t+1}=z^t+\beta(x^{t+1}-z^t)$.
\ENDFOR
	%	\Ensure
\end{algorithmic}
\end{algorithm}

Notice that $\{z^t\}$ in the S-prox-ALM algorithm is an auxiliary sequence, defined as an exponentially weighted average sequence of $\{x^t\}$ .

\subsection{Convergence result}

In the analysis of \cite{zhang2018proximal}, authors assume that $P$ is compact and require both the Slater condition and the strict complementarity condition for the convergence analysis.
 In this paper, we remove these assumptions and prove that we can obtain an $\epsilon$-stationary solution within $\mathcal{O}(1/\epsilon^2)$ iterations using the S-prox-ALM algorithm.

We state our main convergence result, which will be proved in the next two sections.

\begin{theorem}\label{main:NC}
Let $p\ge 3L_f, \rho\ge 0$ and $c<1/(L_f+\rho\smax^2(A)+p)$.
Then there exist $\alpha', \beta'>0$ (depending  on $L_f, A, G, p,\rho, c$ only) such that for all $\alpha< \alpha'$, $\beta< \beta'$ the following results hold:
\begin{enumerate}
\item Every limit point of $\{x^t, y^t\}$ generated by Algorithm \ref{Alg:SProxALM} is a KKT point of \eqref{P1};
\item There exists a constant $B>0$ only depending on $p, \rho, c, \alpha, \beta, L_f, A, G, \underline{f}, f(x^0)$, such that for any $\epsilon>0$, we can find an $\epsilon$-stationary solution within $B/\epsilon^2$ iterations.  In other words, for any $t>0$, we can find an$s\in \{0, 1, \cdots, t-1\}$ such that  $(x^{s+1}, y^{s+1})$ is a $\sqrt{B}/\sqrt{t}$-stationary solution of \eqref{P1}.
\end{enumerate}
\end{theorem}

Note that if we choose $p=3L_f$ and $\rho=L_f$, then the constant $B$ depends on $L_f, A, G, \underline{f}, f(x^0)$ only.
The proof of Theorem~\ref{main:NC} and the explicit representation of $\alpha', \beta'$ will be given in the next section.

\section{Convergence Analysis}

In this section, we use a ``proximal-primal-dual'' framework and a novel global dual error bound to analyze the convergence of the Algorithm~\ref{Alg:SProxALM}.

\subsection{A potential function}
We define
\begin{eqnarray}
x(y, z)&=&\arg\min_{x\in P}K(x, z;y)\nonumber\\
\bar{x}^*(z)&=& \arg\min_{x\in P, Ax=b}\left\{f(x)+\frac{p}{2}\|x-z\|^2\right\},\nonumber\\
P(z)&=&\min_{x\in P, Ax=b}\{f(x)+\frac{p}{2}\|x-z\|^2\}\label{eq:Pdef}\\
d(y, z)&=&\min_{x\in P}K(x, z;y).\label{eq:dK}
\end{eqnarray}

The proof of Theorem~\ref{main:NC} relies on the following potential function:
$$\phi^t=\phi(x^t, y^t, z^t)=K(x^t, z^t;y^t)-2d(y^t, z^t)+2P(z^t).$$
We will prove that $\phi^t$ decreases sufficiently after each iteration of Algorithm~\ref{Alg:SProxALM}, provided that $c, \alpha, \beta$ are chosen sufficiently small.
We start with the following basic descent estimate.

\begin{lemma}\label{four terms}
Let us choose $p,\ c,\ \alpha,\ \beta$ satisfying
\[
p\ge 3L_f,\quad c<1/(L_f+\rho\smax^2(A)+p),\quad \alpha<\frac{c(p-L_f)^2}{4\smax^2(A)},\quad \beta<1/30.
\]

Then for any $t>0$, we have
\begin{equation*}
\phi^t-\phi^{t+1}\ge\frac{1}{4c}\|x^t-x^{t+1}\|^2+\alpha\|Ax(y^{t+1}, z^t)-b\|^2+\frac{p}{3\beta}\|z^t-z^{t+1}\|^2-6p\beta\|x(y^{t+1}, z^t)-\bar{x}^*(z^t)\|^2. \label{key1}
\end{equation*}
\end{lemma}

\noindent{\bf Remark.} It is just the inequality (3.24)  of \cite{zhang2018proximal} and the proof can be seen in \cite{zhang2018proximal}. Though in \cite{zhang2018proximal}, authors assume that $P$ is a bounded box, it is not hard to check this inequality holds for any convex, closed set $P$. For completeness, we will give the proof in appendix \ref{Appendix:A}.
%[[[give the equation number!]]]
%[[Please prove Lemma~\ref{four terms} in the appendix]].

\subsection{The global error bound}

Notice that in Lemma~\ref{four terms}, there is a negative term $\|x(y^{t+1}, z^t)-\bar{x}^*(z^t)\|^2$.  In order to prove that $\phi^t$ is decreasing, we need to bound the term $\|x(y^{t+1}, z^t)-\bar{x}^*(z^t)\|^2$. Notice that if $\|Ax(y^{t+1}, z^t)-b\|=0$, we have $x(y^{t+1}, z^t)=\bar{x}^*(z^t)$. Therefore, it is natural to consider whether we can use $\|Ax(y^{t+1}, z^t)-b\|$ to bound $\|x(y^{t+1}, z^t)-\bar{x}^*(z^t)\|$. Fortunately, it is indeed true. We have the following global dual error bound:

\begin{lemma}\label{dual error bound}
If $p>L_f$, then we have
$$
\|x(y, z)-\bar{x}^*(z)\|<\bar{\sigma}_5\|Ax(y, z)-b\|,\quad \mbox{for any y, z,}
$$
where $\bar{\sigma}_5>0$ depends only on the constants $L=(L_f+\rho\smax^2(A)+p)$, $\gamma=-L_f+p$ and the matrices $A$, $G$.
\end{lemma}

Note that compared to Lemma 3.6 in \cite{zhang2018proximal}, Lemma~\ref{dual error bound} holds globally. Hence, this result is stronger. The constant $\bar{\sigma}_5$ is computable and only depends on the parameters $p, \rho$ and $L_f, A, G$. The explicit estimate of  $\bar{\sigma}_5$ along with the proof of Lemma~\ref{dual error bound} will be given in the next section.

\subsection{Convergence proof}

Equipped with the global dual error bound Lemma \ref{dual error bound}, we prove that the potential function $\phi^t$ decreases sufficiently after each iteration, provided that $\beta$ is chosen to be sufficiently small.
\begin{lemma}\label{sufficient-decrease}
Suppose $p, c, \alpha, \beta$ are chosen to satisfy the conditions in Lemma~\ref{four terms} and we further let
$$\beta<\frac{\alpha}{12p\bar{\sigma}_5^2}.$$
Then for any $t$, we have
\begin{equation}\label{sufficient-descent}
\phi^t-\phi^{t+1}\ge\frac{1}{4c}\|x^t-x^{t+1}\|^2+\frac{\alpha}{2}\|Ax(y^{t+1}, z^t)-b\|^2+\frac{p}{3\beta}\|z^t-z^{t+1}\|^2.
\end{equation}
\end{lemma}

\begin{proof}
By Lemma~ \ref{dual error bound}, we have
\begin{eqnarray*}
6p\beta\|x(y^{t+1}, z^t)-\bar{x}^*(z^t)\|^2
&\le&6p\cdot\frac{\alpha}{12p\bar{\sigma}_5^2}\cdot\bar{\sigma}_5^2\|Ax(y^{t+1}, z^t)-b\|^2\\
&=&\frac{\alpha}{2}\|Ax(y^{t+1}, z^t)-b\|^2.
\end{eqnarray*}
Hence, substituting the above inequality to  Lemma \ref{four terms}  we have
\begin{equation}\label{sufficient-descent1}
\phi^t-\phi^{t+1}\ge\frac{1}{4c}\|x^t-x^{t+1}\|^2+\frac{\alpha}{2}\|Ax(y^{t+1}, z^t)-b\|^2+\frac{p}{3\beta}\|z^t-z^{t+1}\|^2.
\end{equation}
\end{proof}

Next we prove that $\phi^t$ is bounded from below.
\begin{lemma}\label{lower-bounded}
For any $t>0$,
$$\phi^t\ge \underline{f}.$$
\end{lemma}
\begin{proof}
We have
\begin{eqnarray*}
\phi^t
&=&P(z^t)+(K(x^t, z^t;y^t)-d(y^t, z^t))+(P(z^t)-d(y^t, z^t))\\
&\ge&P(z^t)\\
&\ge&\underline{f},
\end{eqnarray*}
where the second step follows from the definition of $d(y^t,z^t)$ (cf.\ \eqref{eq:dK}), and the weak duality $P(z^t)\ge d(y^t,z^t)$, while the last step is due to the boundedness $f$ over the feasible set $\{x\in P\mid Ax=b\}$ (see Assumption~\ref{ass} and \eqref{eq:Pdef}).
\end{proof}

Now we can prove the main theorem.

\medskip
\noindent \textbf{Proof of Theorem \ref{main:NC}: }
We first prove the first part of the main theorem.
For $x, y, z$, we let $F$ be a map such that
$$F(x, y, z)=(x^+, y^+, z^+)$$
is the next iteration point of Algorithm~\ref{Alg:SProxALM}.
It is straightforward to check that the map $F$ is continuous and if $(x,y,z)$ is a fixed point of $F$
$$
F(x,y,z)=(x,y,z),
$$
then $(x,y)\in U^*$ is a pair of primal-dual stationary solution of \eqref{P1}. Suppose that
$$(x^t, y^t, z^t)\rightarrow (\bar{x}, \bar{y}, \bar{z}) \mbox{ along a subsequence $t\in {\cal T}$}.$$
Notice that by Lemma \ref{sufficient-decrease} and Lemma \ref{lower-bounded}, we have
%\begin{eqnarray*}
\[
\|x^t-x^{t+1}\|\rightarrow 0,\quad \|Ax(y^{t+1}, z^t)-b\|\rightarrow 0,\quad
\|z^t-z^{t+1}\|\rightarrow 0.
\]
This further implies
\begin{equation}\label{eq:xyztozero}
\|(x^{t+1}, y^{t+1}, z^{t+1})-(x^t,y^t,z^t)\|\to 0.
\end{equation}
Therefore, we obtain
\begin{eqnarray*}
\|F(\bar{x}, \bar{y}, \bar{z})-(\bar{x}, \bar{y}, \bar{z})\|&=&\lim_{t\rightarrow \infty,\ t\in {\cal T}}\|(x^t, y^t, z^t)-F(x^t, y^t, z^t)\|\\
&=&\lim_{t\rightarrow \infty,\ t\in {\cal T}}\|(x^{t+1}, y^{t+1}, z^{t+1})-F(x^t, y^t, z^t)\|\\
&=&0,
\end{eqnarray*}
where the first step is due to the continuity of $F$ and the second step follows from \eqref{eq:xyztozero}.
Hence, $(\bar{x}, \bar{y})\in X^*$, that is, each limit point $(\bar{x}, \bar{y})$ is a primal-dual stationary solution of \eqref{P1}.

Next we prove that the iteration complexity is $\mathcal{O}(1/\epsilon^2)$.
We need the following ``primal'' error bound, whose proof was given in Lemma 3.6 of \cite{zhang2018proximal} (see also Lemma \ref{error bound} in Appendix \ref{Appendix:A}).

\begin{lemma}\label{lm:eb2}
For any $t\ge 0$, we have
\begin{equation}\label{erb2}
\frac{c(p-L_f)}{1+c(p-L_f)}\|x^{t+1}-x(y^{t+1}, z^t)\|\le \|x^{t+1}-x^t\|.
\end{equation}
\end{lemma}

For $t>0$, we have $\phi^t\ge \underline{f}$.
It follows that
\[
%\begin{eqnarray*}
\sum_{s=0}^{t-1}(\phi^s-\phi^{s+1})=\phi^0-\phi^t
\le\phi^0-\underline{f}.
\]
%\end{eqnarray*}
Hence, there exists an $s\in\{0, \cdots, t-1\}$ such that
\begin{equation}\label{eq:decrease rate}
\phi^s-\phi^{s+1}\le (\phi^0-\underline{f})/t.
\end{equation}
Let
$$C=(\phi^0-\underline{f})\cdot \max\{4c, 2/\alpha, 3\beta/p\}.$$
Then it follows from Lemma~\ref{sufficient-decrease} and \eqref{eq:decrease rate} that
%\begin{eqnarray*}\label{eq:4.3}
\begin{equation}\label{eq:4.3}
\|x^s-x^{s+1}\|^2<C/t,\quad
\|Ax(y^{s+1}, z^s)-b\|^2<C/t,\quad
\|x^{s+1}-z^s\|^2<C/t.
\end{equation}
%\end{eqnarray*}
According to Algorithm \ref{Alg:SProxALM}, we have
$$x^{s+1}= \arg\min_{x}\left\{\langle \nabla_xK(x^s, z^s; y^{s+1}), x-x^s \rangle+\frac{1}{c}\|x-x^s\|^2+\iota(x)\right\}.$$
The corresponding optimality condition is given by
$$0\in \nabla_xK(x^s, z^s;  y^{s+1})+\frac{2}{c}(x^{s+1}-x^s)+\partial{\iota(x^{s+1})}.$$
Letting
$$v=\nabla_xK(x^{s+1}, z^s; y^{s+1})-\nabla_xK(x^s, z^s; y^{s+1})-\frac{2}{c}(x^{s+1}-x^s)-\rho A^T(Ax^{s+1}-b)-p(x^{s+1}-z^s),$$
we can rewrite the above optimality condition as
\[
v\in \nabla_xK(x^{s+1}, z^s; y^{s+1})-\rho A^T(Ax^{s+1}-b)-p(x^{s+1}-z^s)+\partial{\iota(x^{s+1})}.
\]

Recalling the definition \eqref{eq:K} of $K(x,z;y)$, we have
\begin{equation}\label{eq:gradK}
\nabla_x K(x^{s+1},z^s;y^{s+1})=\nabla f(x^{s+1})+A^Ty^{s+1} +\rho A^T(Ax^{s+1}-b)+p(x^{s+1}-z^s).
\end{equation}
Therefore, the optimality condition can be further simplified as
\begin{equation}
v\in \nabla f(x^{s+1})+A^Ty^{s+1}+\partial{(\iota(x^{s+1}))}.
\end{equation}
We now proceed to estimate the size of $v$.
By using the triangle inequality and then using the inequalities \eqref{eq:4.3} and \eqref{erb2}, we have
\begin{eqnarray}
\|Ax^{s+1}-b\|&\le&\|Ax(y^{s+1}, z^s)-b\|+\|A(x^{s+1}-x(y^{s+1}, z^s))\|\nonumber\\
&\le&\frac{\sqrt{C}}{\sqrt{t}}+\smax(A)\frac{1+c(p-L_f)}{c(p-L_f)} \frac{\sqrt{C}}{\sqrt{t}}\nonumber\\
&=&\frac{\sqrt{B_1C}}{\sqrt{t}},\label{eq:4.4}
\end{eqnarray}
where $B_1=\left(1+\smax(A)\frac{1+c(p-L_f)}{c(p-L_f)}\right)^2$.

Furthermore, similar to \eqref{eq:gradK}, we have
\[
\nabla_x K(x^{s},z^s;y^{s+1})=\nabla f(x^{s})+A^Ty^{s+1} +\rho A^T(Ax^{s}-b)+p(x^{s}-z^s),
\]
which can be combined with \eqref{eq:gradK} to obtain
\[
\nabla_x K(x^{s+1},z^s;y^{s+1})-\nabla_x K(x^{s},z^s;y^{s+1})=\nabla f(x^{s+1})-\nabla f(x^{s})+(\rho A^TA+pI)(x^{s+1}-x^s).
\]
This further implies
\[
\|\nabla_x K(x^{s+1},z^s;y^{s+1})-\nabla_x K(x^{s},z^s;y^{s+1})-\frac{2}{c}(x^{s+1}-x^s)\|\le (L_f+p+\rho\smax^2(A)+2/c)\|x^{s+1}-x^s\|,
\]
where we used the Lipschitz continuity of $\nabla f(x)$ (see Assumption~\ref{ass}).
Then we have
\begin{eqnarray*}
\|v\|&\le& (L_f+p+\rho\smax^2(A)+2/c)\|x^s-x^{s+1}\|+\rho\|Ax^{s+1}-b\|+p\|x^{s+1}-z^s\|\\
&\le& (L_f+p+\rho\smax^2(A)+2/c)\frac{\sqrt{C}}{\sqrt{t}}+\rho\smax(A)\frac{\sqrt{B_1C}}{\sqrt{t}}+p\frac{\sqrt{C}}{\sqrt{t}}\\
&\le&\sqrt{B_2C}/\sqrt{t},
\end{eqnarray*}
where the second inequality follows from inequalities \eqref{eq:4.3}, \eqref{eq:4.4}, and
$$B_2=((L_f+p+\rho\smax(A)^2+2/c)+\rho\smax(A)\sqrt{B_1}+p)^2.$$
Notice that $B_1,\ B_2>0$. Then the result holds for $v$ and $B=C\max\{B_1, B_2\}$ and $(x^{s+1}, y^{s+1})$ is a $\sqrt{B}/\sqrt{t}$-stationary solution.
\endprf

\noindent{\bf Remark.} Note that if we take $p=3L_f$, $\rho=L_f$, then the stepsizes need to satisfy the following:
\begin{enumerate}
\item $c<1/(4L_f+L_f\smax^2(A))$;
\item $\alpha<\frac{cL_f^2}{\smax(A)^2}$;
%[[[there is a factor of c in both the denominator and numerator, something is wrong, Check this and Lemma 3.1!!]]]
\item $\beta<\min\{1/30, \frac{\alpha}{12p\bar{\sigma}_5^2}\}$.
\end{enumerate}
Since $\bar{\sigma}_5$ is computable and only depends on $p=3L_f, \rho=L_f, A, G, L_f$, and the stepsizes $c, \alpha, \beta$ only depend on $L_f, A, G$, it follows that the stepsizes are computable and only depend on $A, G, L_f$.
Moreover, according to the proof of the second part of the main theorem, the constant $B$ in the iteration complexity bound only depends  on $\phi^0-\underline{f}, L_f, A, G$.
%%%
\section{A Global Error Bound}

Notice that the key in the proof of the iteration complexity analysis (Theorem~\ref{main:NC}) is the global dual error bound in Lemma~\ref{dual error bound} for the following strongly convex problem (assuming $p>L_f$):
\begin{equation}
\label{P:gz}
\begin{array}{ll}
\mbox{minimize}& f(x)+\frac{\rho}{2}\|Ax-b\|^2+\frac{p}{2}\|x-z\|^2\\ [5pt]
\mbox{subject to} & Ax=b,\ x\in P.
\end{array}
\end{equation}
In this section, we use a decomposition technique to prove this dual error bound for a general strongly convex problem:
\begin{equation}
\label{P:g}
\begin{array}{ll}
\mbox{minimize}& g(x)\\ [5pt]
\mbox{subject to} & Ax=b,\ x\in P,
\end{array}
\end{equation}
where $g$ is a strongly convex function with modulus $\gamma$ and $g$ is smooth with a Lipschitz-continuous gradient of constant $L$.
We define
\begin{eqnarray}
L(x; y)&=&g(x)+y^T(Ax-b),\nonumber\\
x(y)&=&\arg\min_{x\in P}L(x; y);\label{AL-min}\\
d(y)&=&\min_{x\in P}L(x; y),\nonumber\\
x^*&=&\arg\min_{x\in P, Ax=b}g(x),\nonumber\\
p^*&=&g(x^*).\nonumber
\end{eqnarray}
We will prove the following general theorem in this section.
\begin{theorem}\label{eb-linear}
For any $\tilde{y}\in \mathbb{R}^m$, we have
$$\|x(\tilde{y})-x^*\|\le \sigma_5\|Ax(\tilde{y})-b\|,$$
where $\sigma_5>0$ only depends on $L,\ \gamma,\ A,\ G$.
\end{theorem}

%\subsection{Lemma \ref{dual error bound} is a direct corollary of Theorem \ref{eb-linear}}
Before proving Theorem \ref{eb-linear}, we remark that Lemma \ref{dual error bound} is a direct corollary of Theorem \ref{eb-linear}.
In fact, for any $z\in \mathbb{R}^n$, let
$$g(x)=f(x)+\frac{\rho}{2}\|Ax-b\|^2+\frac{p}{2}\|x-z\|^2.$$
Then $g$ is strongly convex with modulus $\gamma=-L_f+p>0$ and smooth with a gradient that is Lipschitz continuous with constant $L=L_f+\rho\smax^2(A)+p$.
Then Lemma~\ref{dual error bound} holds with a constant $\bar{\sigma}_5$ that depends on $p,\ \rho,\ L_f$, and matrices $A,\ G$.

%Finally, we prove Theorem \ref{eb-linear}.
We first sketch the main ideas of the proof of Theorem \ref{eb-linear}.
Let
\begin{eqnarray}
\tilde{r}&=&Ax(\tilde{y})-b, \label{residual} \\
x^*(r)&=&\arg\min_{x: Ax-b=r, x\in P}g(x), \label{proximal with residual}
\end{eqnarray}
where \eqref{proximal with residual} is feasible if and only if $r\in AP-b$.
%[[Should consider the domain of $x^*(\cdot)$ and prove its closeness]]
%\begin{lemma}
%\end{lemma}

To proceed, we write down the KKT conditions of problems \eqref{AL-min} and \eqref{proximal with residual}.
The  KKT condition of problem \eqref{AL-min} are
\begin{eqnarray}\label{KKTforAL}
\nabla g(x(y)) +A^Ty+G^T\mu(y) &=&0, \\
\mu_i(y) &\ge& 0,\ \mbox{for all } i\nonumber\\
x(y)&\in& P,\nonumber\\
\mu_i(y)(Gx(y)-h)_i&=&0,\ \mbox{for all } i.\nonumber
\end{eqnarray}
When $r\in AP-b$, the KKT conditions for \eqref{proximal with residual} are
\begin{eqnarray}\label{KKTforP2}
\nabla g(x^*(r))+A^Ty+G^T\mu&=&0, \\
Ax^*(r)-b&=&r,\nonumber\\
Gx^*(r)&\preceq& h,\nonumber\\
\mu&\succeq& 0,\nonumber\\
\mu_i(Gx^*(r)-h)_i&=&0,\quad i\in[l]\nonumber
\end{eqnarray}
The following lemma shows that problem \eqref{proximal with residual} is always feasible for $r=Ax(y)-b$ with arbitrary $y$.

\begin{lemma}
For any $y$, problem \eqref{proximal with residual} is always feasible with $r=Ax(y)-b$.
Moreover, we have
$$x(y)=x^*( {r}),\ x^*( 0)=x^*$$
for $r=Ax(y)-b$.
\end{lemma}

\begin{proof}
Notice that if we add the condition $Ax(y)-b=r$ to the KKT system \eqref{KKTforAL}, we just attain the KKT system \eqref{KKTforP2}.

It means that $(x(y), y)$ are the solution to the KKT system of problem \eqref{proximal with residual}. Due to the strong convexity, the primal solution of the problem \eqref{proximal with residual} is unique. Hence, we have $x(y)=x^*( {r})$.
The claim $x^*(0)=x^*$ follows directly from the definition.
\end{proof}

\begin{lemma}
Let $\tilde{r}=Ax(\tilde{y})-b\in AP-b$.
Then for any $s\in[0, 1]$ and $r=s\tilde{r}$, the problem \eqref{proximal with residual} is feasible, i.e., $s\tilde{r}\in AP-b$.
Hence, $x^*( r)$ is well-defined for $r=s\tilde{r}, s\in[0, 1]$.
\end{lemma}
\begin{proof}
In fact, since $x^*( \tilde{r})$ satisfies
$Ax^*( \tilde{r})-b=\tilde{r}, x^*( \tilde{r})\in P$ and $x^*$ satisfies $Ax^*-b=0, x^*\in P$, we have
\begin{eqnarray}
A(sx^*( \tilde{r})+(1-s)x^*)-b
&=&s(Ax^*( \tilde{r})-b)+(1-s)(Ax^*-b)\nonumber\\
&=&s\tilde{r}\label{eq2}\\
&=&r\nonumber,
\end{eqnarray}
where \eqref{eq2} is because $Ax^*=b$. Moreover, $sx^*( \tilde{r})+(1-s)x^*\in P$ since $P$ is convex. Hence $sx^*( \tilde{r})+(1-s)x^*$ is feasible to the problem \eqref{proximal with residual} with $r=s\tilde{r}$ and $r\in AP-b$.
Then since $g(x)$ is strongly convex, $\arg\min_{x\in P, Ax=b}g(x)$ must have a solution. Hence, $x^*(r)$ is well-defined.
\end{proof}
Let $Y^*(r)$ be the dual solution set of the problem \eqref{proximal with residual}.
For arbitrary $\tilde{y}\in \mathbb{R}^m$, the idea of the proof is to find finite scalar $0=\eta_0<\eta_1<\eta_2<\cdots<\eta_R=1$ and $r^i=\eta_i\tilde{r}$, such that
$$\|x^*( r^i)-x^*( r^{i-1})\|\le \sigma_5\|r^i-r^{i-1}\|.$$
Then summing up this inequality from $1$ to $R$ and using the fact that
\begin{eqnarray}\label{final}
\sum_{i=1}^R\|(r^i-r^{i-1})\|&=&\sum_{i=1}^R(\eta_i-\eta_{i-1})\|\tilde{r}\|\\
&=&\eta_R\|\tilde{r}\|\nonumber\\
&=&\|A(x(\tilde{y})-x^*)\|\nonumber\\
&=&\|Ax(\tilde{y})-b\|\nonumber
\end{eqnarray}
and
\[
\|x(\tilde{y})-x^*\|=\|x^*(r^{R})-x^*(r^0)\|\le \sum_{i=1}^R\|x^*( r^i)-x^*( r^{i-1})\|,
\]
we attain the desired result.
Therefore, the decomposition step is the key to the whole proof.
\subsection{The proof of Theorem \ref{eb-linear}}
To prove Theorem \ref{eb-linear}, we need a series of lemmas.
The following is the well known Hoffman bound which estimates the distance from a point to a polyhedral set by the amount of infeasibility.
\begin{proposition}\label{Hoffman}
Let $C_1\in \mathbb{R}^{m_1\times n},\ C_2\in \mathbb{R}^{m_2\times n}$ and $b_1\in \mathbb{R}^{m_1} ,b_2\in \mathbb{R}^{m_2}$, then the distance from a point %[[\bar{x}]]$x\in \mathbb{R}^n$ to the set:
$$S=\{x\mid C_1x\preceq  b_1, C_2x=b_2\}$$
is bounded by
$$\dist(\bar{x}, S)^2\le\theta(\|(C_1\bar{x}-b_1)_+\|^2+\|C_2\bar{x}-b_2\|^2),$$
where $(\cdot)_+$ means the projection to the nonnegative orthant and $\theta$ is a positive constant depending  on $C_1$ and $C_2$ only.
\end{proposition}
Here we call $\theta$ as the Hoffman constant of the linear system.
The Hoffman constant $\theta$ is computable. The proof of Proposition~\ref{Hoffman} and the estimate of $\theta$ can be found in Lemma 3.2.3 of \cite{Pang-complementary}. The following lemma shows the continuity of $x^*(r)$ as a map of $r$.
\begin{lemma}\label{continuity}
$x^*(r)$ is continuous of $r$ for $r\in AP-b$. In other words, for any $r\in AP-b$ and a sequence $\{v^i\}\subseteq AP-b$ with $v^i\rightarrow r$, we have $x^*(v^i)\rightarrow x^*(r)$.
\end{lemma}

\begin{proof}
First, we prove that $x^*(v^i)$ is bounded .
We prove it by contradiction. Assume that $\|x^*(v^i)\|\rightarrow \infty$.
Then by the coerciveness of $g$, we  have  $g(x^*(v^i))\rightarrow \infty$.
On the other hand, let $x^i$ be the projection of $x^*(r)$ to the set $\{x\mid Ax-b=v^i, x\in P\}$.
Then, by Hoffman bound (Proposition~\ref{Hoffman}), we know that $\|x^*(r)-x^i\|<\sqrt{\theta}\|r-v^i\|\rightarrow 0$.
Hence, $g(x^i)\rightarrow g(x^*(r))$, implying $\{g(x^i)\}$ is bounded. Since $g(x^*(v^i))\to\infty$, it follows that for $i$ large enough, $g(x^i)<g(x^*(v^i))$. This is a contradiction to the definition of $x^*(v^i)$.
Therefore, $\{x^*(v^i)\}$ is bounded.
Next, we prove $x^*(v^i)\to x^*(r)$ by contradiction. Assume the contrary that $x^*(v^i)$ is a sequence that does not converge to $x^*(r)$.
Since $\{x^*(v^i)\}$ is bounded, passing to a subsequence if necessary,  we assume that $x^*(v^i)\rightarrow \bar{x}\in P$ and $\bar{x}\ne x^*(r)$.
Since the active set of any $x^*(v^i)$ is a subset of $[l]$, there exists a subsequence of $\{x^*(v^i)\}$ that has a common active set $\mathcal{A}$, i.e.,
we have the KKT conditions:
\begin{eqnarray*}
\nabla g(x^*(v^i))+A^Ty+G^T\mu&=&0\\
Ax^*(v^i)-b&=&v^i\\
(Gx^*(v^i))_j&=&h_j,\ j\in \mathcal{A}\\
Gx^*(v^i)&\preceq &h,\\
\mu_j&=& 0,\ j\in \mathcal{A}\\
\mu_j&\ge&0,\ j\notin \mathcal{A}.
\end{eqnarray*}
Define $T=\{A^Ty+G^T\mu\mid \mu_j\ge0,\ j\in \mathcal{A};\ \mu_j= 0,\ j\notin \mathcal{A}\}$. Note that $T$ is a finitely generated cone, hence it is closed (seeProposition 3.2.1(a) in \cite{bertsekas-cvx}). Taking limit in the KKT system, we have
$$-\nabla g(\bar{x})\in T\ \ \mbox{
and }\ \ A\bar{x}-b=r.$$
Consequently, we see that $\bar{x}$ satisfies the KKT conditions of the problem
$$\min_{Ax-b=r, x\in P}g(x).$$
By strong convexity, the solution of the above problem is unique, hence, $\bar{x}=x^*(r)$, which  is a contradiction. This completes the proof.
\end{proof}

We define a set-valued function (also called multifunction) for $r\in AP-b$:
\begin{equation}
M(r)=\{(y, \mu)\mid (x^*(r), y, \mu)\mbox{ is a primal-dual solution of \eqref{proximal with residual}}\},
\end{equation}
which maps a vector $r$ to a set in $\mathbb{R}^m\times \mathbb{R}^l$.
\begin{definition}
$\mathcal{A}$ is said to be an active set of $r$ if there exists a $w=(y, \mu)\in M(r)$ such that
\begin{eqnarray*}
\nabla g(x^*( r))+A^Ty+G^T\mu&=&0\\
Ax^*( r)-b&=&r\\
(Gx^*(r))-h)_j&=&0,\mbox{ \ }j\in \mathcal{A}\\
(Gx^*( r)-h)_j&\le&0, \mbox{ \ } j\notin \mathcal{A}\\
\mu_j&\ge&0, \mbox{ \ }j\in \mathcal{A}\\
\mu_j&=&0, \mbox{ \ }j\notin \mathcal{A}.
\end{eqnarray*}
%We also say that $\mathcal{A}$ is an active set of $w$.
Moreover, we say $r$ and $r'$ share a common active set $\mathcal{A}$ if $\mathcal{A}$ is an active set of both $r$ and $r'$.
% there exists a set $\mathcal{A}$, $w\in M(r)$ and $w'\in M(r')$ such that
%$\mathcal{A}$ is an active set of $w$ and $w'$.
%In this case, we also say that $w, w'$  share a common active set.
%[[simplify the notations]]
\end{definition}
%\begin{lemma}\label{equivalent}
%$w\in M(r)$ and $w'\in M(r')$ share a common active set if and only if the following two conditions hold:
%\begin{enumerate}
%\item For any $j\in [l]$ with $(Gx^*( r)-h)_j<0$, we always have $\mu_j'=0$.
%\item For any $j$ with $\mu_j'>0$, we always have $(Gx^*( r')-h)_j=0$.
%\end{enumerate}
%\end{lemma}
%#
%The proof is trivial and is omitted.
We will then state some properties of $M(\cdot)$, which will be used to prove the dual error bound.

\begin{lemma}\label{general Hoffman}
Suppose $r,\ r'$ share a common active set $\mathcal{A}$. Then there exist $w\in M(r)$ and $w'\in M(r')$ with the same active set $\mathcal{A}$ such that
\begin{equation}
\dist(M(r'), M(r))+\|x^*(r)-x^*(r')\|\le \|w-w'\|+\|x^*(r)-x^*(r')\|\le \sigma_5\|r-r'\|
\end{equation}
for some $\sigma_5>0$ which depends only on $L,\ \gamma,\ A,\ G$.
\end{lemma}
\begin{proof}
%Without loss of generality, suppose the projection of $w'$ to $M(r)$ is $w=(y, \mu)$.
Let $\mathcal{A}$ be the active set shared by $r$ and $r'$. Then, by definition, there exist $w=(y,\mu)\in M(r)$ and $w'=(y',\mu')\in M(r')$ such that the KKT conditions for $x^*( r)$ and $x^*(r')$ can be written as
\begin{eqnarray*}
A^Ty+G^T\mu&=&-\nabla g(x^*( r))\\
Ax^*( r)-b&=&r\\
\mu_j&=&0,\quad j\notin \mathcal{A}\\
\mu_j&\ge&0,\quad j\in \mathcal{A}\\
(Gx^*(r)-h)_j&=&0,\quad j\in \mathcal{A}\\
Gx^*( r)-h&\preceq& 0.
\end{eqnarray*}
and
\begin{eqnarray*}
A^Ty'+G^T\mu'&=&-\nabla g(x^*( r'))\\
Ax^*( r')-b&=&r'\\
\mu_j'&=&0,\quad j\notin \mathcal{A}\\
\mu_j'&\ge&0,\quad j\in \mathcal{A}\\
(Gx^*(r')-h)_j&=&0,\quad j\in \mathcal{A}\\
Gx^*( r')-h&\preceq & 0.
\end{eqnarray*}
Notice that $w'=(y', \mu')$ satisfies the linear system
\begin{eqnarray*}
A^Ty+G^T\mu&=&-\nabla g(x^*( r))\\
\mu_j&=&0, \quad j\notin \mathcal{A}\\
\mu_j&\ge&0,\quad j\in \mathcal{A}
\end{eqnarray*}
approximately, it follows from the Hoffman bound (Proposition~\ref{Hoffman}) that
\begin{eqnarray}\nonumber
\|w'-w\|^2+\|x^*(r)-x^*(r')\|^2 &\le& \theta\|\nabla g(x^*( r))-\nabla g(x^*( r'))\|^2+\|x^*(r)-x^*(r')\|^2\\
&
%\stackrel{ \mbox{\scriptsize(i)}}
\le&
(\theta L^2+1)\|x^*(r)-x^*(r')\|^2\nonumber\\
&\le &\frac{\theta L^2+1}{\gamma}\langle \nabla g(x^*( r))-\nabla g(x^*( r')), x^*( r)-x^*(r')\rangle,\label{intermediate}
\end{eqnarray}
for some $\theta$ which depends only on $A,\ G,\ \mathcal{A}$, where the second step follows from the Lipschitz continuity of $\nabla g(x)$, and the last step is due to the strong convexity of $g$.
%because for $x, x'$, we have
%\begin{eqnarray*}
%\|\nabla g(x)-\nabla g(x')\|^2&\le& L^2\|x-x'\|^2\\
%&\le& \frac{L^2}{\mu_g}\langle \nabla g(x)-\nabla g(x'), x-x'\rangle.
%\end{eqnarray*}

%Here $\theta$ only depends on $A, G, \mathcal{A}$, which is coming from the Hoffman bound in Proposition \ref{Hoffman}.
Notice that the KKT conditions imply
$$\nabla g(x^*( r))=-A^Ty-G^T\mu, \ \  \nabla g(x^*( r'))=-A^Ty'-G^T\mu',\ \ A(x^*( r)-x^*( r'))=r-r'$$
and
\begin{equation*}
 \mu'_{\bar{\mathcal{A}}}=\mu_{\bar{\mathcal{A}}}=0, \quad (Gx^*(r')-h)_{\mathcal{A}}=(Gx^*(r')-h)_{\mathcal{A}}=0.
%\langle \mu', Gx^*(r)-h\rangle\le 0, \ \
%\langle \mu, Gx^*(r)-h\rangle= 0,\ \
%\langle \mu', Gx^*(r')-h\rangle= 0.
\end{equation*}
The latter implies
\[
(\mu'-\mu)_{\bar{\mathcal{A}}}=0,\quad (G(x^*(r)-x^*(r')))_{\mathcal{A}}=0
\]
so that
\[
\langle G^T(\mu'-\mu), x^*(r)-x^*(r')\rangle =\langle (\mu'-\mu)_{\mathcal{A}}, (G(x^*(r)-x^*(r')))_{\mathcal{A}}\rangle +\langle (\mu'-\mu)_{\bar{\mathcal{A}}}, (G(x^*(r)-x^*(r')))_{\bar{\mathcal{A}}}\rangle= 0.
\]
Hence, we have
\begin{eqnarray*}
\langle \nabla g(x^*(r))-\nabla g(x^*(r')), x^*(r)-x^*(r')\rangle&=&=\langle A^T(y'-y)+G^T(\mu'-\mu), x^*(r)-x^*(r')\rangle \\
&=& \langle y-y',A(x^*( r)-x^*( r'))\rangle\\
&=& \langle y-y',r-r'\rangle.
\end{eqnarray*}
Combining this with \eqref{intermediate}, we obtain
\begin{eqnarray*}
\|w'-w\|^2+\|x^*(r)-x^*(r')\|^2&\le& \frac{\theta L^2+1}{\gamma}\|y-y'\|\|r-r'\|\\
&\le&\frac{\theta L^2+1}{\gamma}\|w-w'\|\|r-r'\|\\
&\le &\frac{\theta L^2+1}{\gamma}(\|w-w'\|^2+\|x^*(r)-x^*(r')\|^2)^{1/2}\|r-r'\|,
\end{eqnarray*}
further implying
\[
(\|w'-w\|^2+\|x^*(r)-x^*(r')\|^2)^{1/2}\le \frac{\theta L^2+1}{\gamma}\|r-r'\|.
\]
Finally, we have
\begin{eqnarray*}
\dist(M(r'),M(r))+\|x^*(r)-x^*(r')\|&\le& \|w'-w\|+\|x^*(r)-x^*(r')\|\\
&\le &\sqrt{2}(\|w'-w\|^2+\|x^*(r)-x^*(r')\|^2)^{1/2}\\
&\le& \frac{\sqrt{2}(\theta L^2+1)}{\gamma}\|r-r'\|.
\end{eqnarray*}
Notice that $\theta$  depends only on $A, G, \mathcal{A}$. Since $\mathcal{A}\subseteq[l] $ has only finitely many choices, we can take $\bar{\theta}$ to be the maximum Hoffman constant for all  linear systems (defined by the subset $\mathcal{S}\in \subseteq[l]$) of the form
\begin{eqnarray*}
A^Ty+G^T\mu&=&-\nabla g(x^*( r))\\
\mu_j&=&0, \quad j\notin \mathcal{S}\\
\mu_j&\ge&0,\quad j\in \mathcal{S}.
\end{eqnarray*}
%is $\theta(A, G, \mathcal{S})$, we let
%\begin{equation}\label{bartheta}
%\bar{\theta}=\max_{S\subseteq [l]}\theta(A, G, \mathcal{S}).
%\end{equation}
%Then we let
%$$\sigma_5=\bar{\theta}\frac{L^2}{\gamma}$$
Then we have
$$\dist(M(r'), M(r))+\|x^*(r)-x^*(r')\|\le \|w-w'\|+\|x^*(r)-x^*(r')\| \le \sigma_5\|r-r'\|$$
where
\[
\sigma_5= \frac{\sqrt{2}({\bar{\theta}} L^2+1)}{\gamma}
\]
which depends only on $L,\ \gamma,\ A,\ G$. The estimation of $\bar{\theta}$ will be discussed in Section~\ref{stepsize estimate}.
\end{proof}

Lemma~\ref{general Hoffman} can be considered a kind of Lipschitzian property of the multifunction $M(r)$. The following lemma states some additional continuity property of $M(\cdot)$.
%[[I add (c)]]
\begin{lemma}\label{properties}
The multifunction $M(\cdot)$ has the following properties:
\begin{enumerate}
\item [{\rm (a)}] $M(r)$ is a closed set for any $r\in AP-b$.
\item [{\rm (b)}] Suppose that $\{v^i\}\subseteq AP-b$ is a sequence converging to $r\in AP-b$ and for any $i$, $v^i$ shares a common active set with $r'$. Then $r$ shares a common active set with $r'$.
\item [{\rm (c)}] If $r\in AP-b$ and $\{v^i\}\subseteq AP-b$ with $v^i\rightarrow r$. If $\mathcal{A}$ is an active set of any $v^i$, then $\mathcal{A}$ is also an active set of $r$.
\end{enumerate}
\end{lemma}
\begin{proof}
The proof of first claim is straightforward by checking the KKT conditions and using the fact a finitely generated cone is closed (see Proposition 3.2.1(a) in \cite{bertsekas-cvx}).
We now prove the second claim.
Since the choice for $\mathcal{A}$ is finite, there exists  a sub-sequence of $\{v^i\}$ sharing a common active set $\mathcal{A}$ with $r'$.
To simplify notation, let us denote this subsequence still by $\{v^i\}$. % shares a common active set $\mathcal{A}$ with $r'$.
Let
$$T=\{A^Ty+G^T\mu\mid \mu_j\ge0, j\in \mathcal{A}; \mu_j= 0, j\notin \mathcal{A}\}.$$
Then by the definition of $\mathcal{A}$ and $M(\cdot)$, we have
$$-\nabla g(x^*(v^i))\subseteq T$$
for any $i$.
Notice that $T$ is a finitely generated cone, hence $T$ is closed (see Proposition 3.2.1(a) in \cite{bertsekas-cvx}).
On the other hand, by Lemma~\ref{continuity}, we have
$$x^*(r)=\lim_{i\rightarrow \infty}x^*(v^i).$$
It follows from the closeness of $T$  and the continuity of $x^*(\cdot)$ (Lemma \ref{continuity}) that
$$-\nabla g(x^*(r))=-\lim_{i\rightarrow \infty}\nabla g(x^*(v^i))\in T.$$
By definition of $T$, there exists a $(y, \mu)$ satisfying the following system
\begin{eqnarray*}
\nabla g(x^*( r))+A^Ty+G^T\mu&=&0,\\
Ax^*( r)-b&=&r,\\
(Gx^*( r))_j-h&=&0, \quad j\in \mathcal{A}\\
(Gx^*( r)-h)_j&\le&0, \quad j\notin \mathcal{A}\\
\mu_j&\ge&0, \quad j\in \mathcal{A}\\
\mu_j&=&0, \quad j\notin \mathcal{A},
\end{eqnarray*}
which proves the desired result.
The third  claim is just a direct corollary of the second claim. In fact, letting $r'=v^1$ in the second claim, we attain that $\mathcal{A}$ is a common active set of $v^1$ and $r$. This implies that $\mathcal{A}$ is an active set of $r$.
\end{proof}
The following lemma shows that for a fixed $r\in AP-b$, the set of all vectors in $AP-b$ sharing a common active set with $r$ is a closed set.
%[[I add a lemma]]
\begin{lemma}\label{closed-set}
For any $r\in AP-b$, the set $Q(r)\{r'\in AP-b\mid \mbox{r\ and\ r'\ share\ a\ common\ active\ set}\}$ is closed.
\end{lemma}
\begin{proof}
First, since the image of a polyhedral  under a linear transformation is still a polyhedral (see exercise 3.9(b) in \cite{bertsekas-cvx}), the set $AP-b$ is closed.
Suppose that a sequence $\{v^i\}\subseteq AP-b$ satisfies that $v^i\rightarrow v$ and $v^i$ shares a common active set with $r$. Since $AP-b$ is closed ,
it follows that $v\in AP-b$. Then by Part (b) in Lemma \ref{properties}, $v$ shares a common active set with $r$. Hence, $Q(r)$ is closed.
\end{proof}
We are now ready to prove the following proposition which is the key for the global error bound. It enables us to decompose the line segment from $0$ to $\tilde{r}$ to finitely many local pieces and the error bound holds in each piece.
\begin{proposition}\label{dual-pieces}
For any $\tilde{r}$ be given by \eqref{residual}, there exists a sequence $0=\eta_0<\eta_1<\cdots<\eta_R=1$ such that for any $i\in \{0, 1,\cdots, R-1\}$, $\eta_i\tilde{r}$ and $\eta_{i+1}\tilde{r}$ share a common active set.
\end{proposition}

\begin{proof}
We set $\eta_0=0$ and then construct $\eta_i$ recursively.
If we have already defined $\eta_{i-1}<1$, we define $\eta_i$ to be the largest $\eta\in [\eta_{i-1},1]$ such that $\eta\tilde{r}$ shares a common active set with $\eta_{i-1}\tilde{r}$.
We first need to verify that this ``largest'' $\eta$ exists.
In fact, by Lemma \ref{closed-set}, the set
$$Q(\eta_{i-1}\tilde{r}):=\{r \mid {r}\mbox{ shares\ a\ common\ active\ set\ with }\eta_{i-1}\tilde{r}\}$$
is closed. Hence, $Q(\eta_{i-1}\tilde{r})\cap\{r=\eta \tilde{r}\mid \eta\in[\eta_{i-1}, 1]\}$ is compact. Therefore,
the ``largest'' $\eta$ exists.

Next, we need to prove that if $\eta_{i-1}<1$, we have $\eta_i>\eta_{i-1}$.
Equivalently, we need to show that there exists an $\eta\in(\eta_{i-1}, 1]$ such that $\eta_{i-1}\tilde{r}$ shares a common active set with $\eta\tilde{r}$.
In fact, let $\{\lambda^j\}\subseteq (\eta_{i-1}, 1]$ be an arbitrary  sequence satisfying $\lambda^j\rightarrow \eta_{i-1}$ and $\lambda^j>\eta_{i-1}$.
Then there exists a subsequence $\lambda^{j_k}$ such that $\lambda^{j_k}\tilde{r}$ shares a common active set $\mathcal{A}$ with $\lambda^{j_1}\tilde{r}$.
By Lemma \ref{properties}(c), $\eta_{i-1}\tilde{r}$ shares a common active set with $\lambda^{j_1}\tilde{r}$. Hence,
$\mathcal{A}$ is also an active set of $\eta_{i-1}\tilde{r}$. Therefore there exists an $\eta\in(\eta_{i-1}, 1]$ such that $\eta\tilde{r}$ sharing a common active set with $\eta_{i-1}\tilde{r}$(here $\eta=\lambda^{j_1}$).
Thus, $\eta_i>\eta_{i-1}$.

Finally, We need to prove that there exists an $R<\infty$ such that $\eta_R=1$. Let $r^i=\eta_i\tilde{r}$.
We prove by contradiction. Assume that $\eta_i<1$ for any positive integer $i$. Since the choice for an active set is finite, there exist $\eta_{i_1}<\eta_{i_2}<\eta_{i_3}$ such that $r^{i_1}, r^{i_2}, r^{i_3}$ share a common active set $\mathcal{A}$. Clearly, $\eta_{i_3}>\eta_{i_2}\ge \eta_{i_1+1}$. However, if $r^{i_1}, r^{i_3}$ share a common active set, then by the definition of $\eta_{i_1+1}$ we should have $\eta_{i_1+1}\ge \eta_{i_3}$. This is a contradiction.
\end{proof}

Now we are ready to finish the proof of Theorem \ref{eb-linear}.

\noindent
{\bf Proof of Theorem \ref{eb-linear}.}

By Proposition \ref{dual-pieces}, we can define $r^i=\eta_i\tilde{r}$, $i=0,1,...,R-1, R$, so that
$r^{i+1}$ and $r^i$ share a common active set for $i=0, \cdots, R-1$. It follows from Lemma~\ref{general Hoffman} that, for $1\le i\le R-1$, we have
\begin{eqnarray}
\|x^*( r^i)-x^*( r^{i-1})\|%&\le & \dist(M(r^i), M(r^{i+1}))\\
&\le& \sigma_5\|r^i-r^{i-1}\|. \label{piecebound}
\end{eqnarray}
Notice that by the definition of $r^i$, we have
\begin{eqnarray}
\sum_{i=1}^R\|(r^i-r^{i-1})\|&=&\sum_{i=1}^R(\eta_i-\eta_{i-1})\|\tilde{r}\|\\
&=&\eta_R\|\tilde{r}\|\\
&=&\|A(x(\tilde{y})-x^*)\|\\
&=&\|Ax(\tilde{y})-b\| \label{sum1}
\end{eqnarray}
and
\begin{equation}\label{sum2}
\|x(\tilde{y})-x^*\|=\|x^*(r^{R})-x^*(r^0)\|\le \sum_{i=1}^R\|x^*( r^i)-x^*( r^{i-1})\|,
\end{equation}
where we use the triangle inequality in the last inequality.
Summing inequality \eqref{piecebound} from $1$ to $R$, and using \eqref{sum1} and \eqref{sum2}, we finish the proof.
\endprf

\subsection{Stepsize estimate}\label{stepsize estimate}

According to \cite{Pang-complementary}, we can compute Hoffman's error bound constant ${\theta}$ as follows:
Let
$$\mathbf{M}=\left[\begin{array}{cc}
A^T&G^T\\
0&I
\end{array}\right].$$
Then
$$\bar{\theta}=\max_{\mathbf{\bar{M}}\in \mathcal{B}(\mathbf{M})}\sigma_{\max}^2(\bar{M})/\sigma_{\min}^4(\bar{M}),$$
where $\mathcal{B}(\mathbf{M})$ is the set of all sub-matrices of $\mathbf{M}$ with full row rank.
Note that
\begin{equation*}
{\sigma}_5=\frac{\sqrt{2}(\bar{\theta} L^2+1)}{\gamma}=\frac{\sqrt{2}(\bar{\theta} (L_f+\rho\smax(A)+p)^2+1)}{-L_f+p}.
\end{equation*}
Thus, the stepsizes satisfying
\[
c<1/(4L_f+L_f\smax^2(A)),\quad \alpha<\frac{cL_f^2}{\smax^2(A)},\quad \beta<\min\{1/30, \frac{\alpha}{12p{\sigma}_5^2}\}
\]
are computable and, according to Theorem~\ref{main:NC},  will ensure the convergence of the S-prox-ALM algorithm.
\appendix
\appendixpage
\addappheadtotoc
\section{Proof of Lemma \ref{four terms}}\label{Appendix:A}
In this section, we prove Lemma \ref{four terms}.
To do this , we need a series of lemmas which are proved in \cite{zhang2018proximal}.
The first lemma states some primal error bounds.
\begin{lemma}[Error Bounds]\label{error bound}
Suppose $p>L_f$, $\rho>0$ are fixed. Then there exist positive constants $\sigma_1,...,\sigma_4 >0$ $($independent of $y$ and $z$ such that the following error bounds hold:
\begin{eqnarray}
\|x^{t+1}-x^t\|&\ge& \sigma_1\|x^t-x(y^{t+1}, z^t)\|,\label{eb1}\\
\|x^{t+1}-x^t\|&\ge& \sigma_2\|x^{t+1}-x(y^{t+1}, z^t)\|, \label{eb2} \\
\|y-y'\|&\ge&\sigma_3\|x( y,z)-x(y',z)\|, \label{eb6}\\
  \|z^t-z^{t+1}\| &\ge &{\sigma}_4 \|\bar{x}^*(z^t)-\bar{x}^*(z^{t+1})\|,\label{eb4} \\
 \|z^t-z^{t+1}\|&\ge &  \sigma_4\|x(y^{t+1}, z^t)-x(y^{t+1}, z^{t+1})\|,\label{eb3}
\end{eqnarray}
for any $y, y'$.
%[[]]\tilde{\sigma}_4
where $\sigma_1=c\gamma_K=c(p-L_f)$, $\sigma_2={\sigma_1}/({1+\sigma_1})$, $\sigma_3=\gamma_K/\smax(A)=(-L_f+p)/\smax(A)$, ${\sigma}_4=\gamma_K/p=(-L_f+p)/p$.
\end{lemma}

As in \cite{zhang2018proximal}, we have the following three descent lemmas.
\begin{lemma}[Primal Descent] \label{primal}
For any t, if $c<1/(L_f+\rho\smax^2(A)+p),$
$$K(x^t, z^t; y^t)-K(x^{t+1}, z^{t+1}; y^{t+1})\ge \frac{1}{2c}\|x^t-x^{t+1}\|^2+\frac{p}{2\beta}\|z^t-z^{t+1}\|^2-\alpha\|Ax^t-b\|^2.$$
\end{lemma}
\begin{lemma}[Dual Ascent]\label{dual ascent}
For any $t$, we have
\begin{equation}\nonumber
d(y^{t+1}, z^{t+1})-d(y^t, z^t)
\ge  \alpha(Ax^t-b)^T(Ax(y^{t+1}, z^t)-b)
+\frac{p}{2}(z^{t+1}-z^t)^T(z^{t+1}+z^t-2x(y^{t+1}, z^{t+1}))
\end{equation}
\end{lemma}
\begin{lemma}[Proximal Descent]\label{proximal-descent}
For any $t\ge 0$, there holds
\begin{equation}\label{eq:prox-descent}
P(z^{t+1})-P(z^{t})\le p(z^{t+1}-z^t)^T(z^t-\bar{x}^*(z^t))+\frac{p}{2{\sigma}_4}\|z^t-z^{t+1}\|^2.
\end{equation}
\end{lemma}

\noindent{\bf Remark.} Note that in \cite{zhang2018proximal}, $P$ represents a bounded box. However, as  stated in the remark after Lemma 3.6 in \cite{zhang2018proximal}, the lemmas above hold for any closed convex set $P$.

We now use the above lemmas to prove Lemma \ref{four terms}.

\noindent
{\bf Proof of Lemma \ref{four terms}}
Recall the definition of potential function
$$\phi^t=K(x^t, z^t;y^t)-2d(y^t, z^t)+2P(z^t).$$
Combining the above lemmas, we have
\begin{eqnarray}
&&\!\!\!\!\!\!\!\!\!\!\!\!\phi^t-\phi^{t+1}\nonumber\\
&\ge & \left(\frac{1}{2c}\|x^{t+1}-x^t\|^2-\alpha\|Ax^t-b\|^2+\frac{p}{2\beta}\|z^t-z^{t+1}\|^2\right)\nonumber\\
&&+2 \left(\alpha(Ax^t-b)^T(Ax(y^{t+1}, z^t)-b)+p(z^{t+1}-z^t)^T(z^{t+1}+z^t-2x(y^{t+1}, z^{t+1})) \right)\nonumber\\
&&+2 \left(p(z^{t+1}-z^t)^T(\bar{x}^*(z^t)-z^t)-\frac{p}{2{\sigma}_4}\|z^{t+1}-z^t\|^2 \right)\nonumber\\
&=& \left(\frac{1}{2c}\|x^{t+1}-x^t\|^2-\alpha\|Ax^t-b\|^2+\frac{p}{2\beta}\|z^t-z^{t+1}\|^2\right)
+2\alpha(Ax^t-b)^T(Ax(y^{t+1}, z^t)-b)\nonumber\\
&&+p(z^{t+1}-z^t)^T\left((z^{t+1}-z^t)-2(x(y^{t+1}, z^{t+1})-\bar{x}^*(z^t))\right)
-\frac{p}{2{\sigma}_4}\|z^t-z^{t+1}\|^2\nonumber\\
&=& {\color{black}\left(\frac{1}{2c}\|x^{t+1}-x^t\|^2-\alpha\|Ax^t-b\|^2+\frac{p}{2\beta}\|z^t-z^{t+1}\|^2\right)
+2\alpha(Ax^t-b)^T(Ax(y^{t+1}, z^t)-b)}\nonumber\\
&&+p(z^{t+1}-z^t)^T\left((z^{t+1}-z^t)-2(x(y^{t+1}, z^{t+1})-x(y^{t+1}, z^t))+2(x(y^{t+1}, z^t)-\bar{x}^*(z^t))\right)\nonumber\\
&&-\frac{p}{2{\sigma}_4}\|z^t-z^{t+1}\|^2.\label{first descent}
\end{eqnarray}
Let $\zeta$ be an arbitrary positive scalar, and by the fact that
$$\|(z^{t+1}-z^t)/\zeta+\zeta(x(y^{t+1}, z^t)-\bar{x}^*(z^t))\|^2\ge 0,$$
we have
$$2(z^{t+1}-z^t)^T(x(y^{t+1}, z^t)-\bar{x}^*(z^t))\ge -\|z^t-z^{t+1}\|^2/\zeta-\zeta\|x(y^{t+1}, z^{t})-\bar{x}^*(z^t)\|^2.$$
Using Cauchy-Schwarz inequality and the error bound \eqref{eb3} in Lemma~\ref{error bound}, we have
%[[]]missing 2
\begin{eqnarray*}
-2(z^{t+1}-z^t)^T(x(y^{t+1}, z^{t+1})-x(y^{t+1}, z^t))&\ge& -\|z^t-z^{t+1}\|\|x(y^{t+1}, z^{t+1})-x(y^{t+1}, z^t)\|\\
&\ge& -\frac{1}{\sigma_4}\|z^t-z^{t+1}\|^2.
\end{eqnarray*}
Substituting these two inequalities into \eqref{first descent}, we have
\begin{eqnarray*}
&&\phi^t-\phi^{t+1}\nonumber\\
&\ge & \frac{1}{2c}\|x^{t+1}-x^t\|^2-\left(\alpha\|Ax^t-b\|^2-2\alpha(Ax^t-b)^T(Ax(y^{t+1}, z^t)-b)+\alpha\|Ax(y^{t+1}, z^t)-b\|^2\right)\nonumber\\
&&+\alpha\|Ax(y^{t+1}, z^t)-b\|^2%\nonumber\\
%&&
+\left(\frac{p}{2\beta}+p-\frac{p}{\sigma_4}-\frac{p}{\zeta}-\frac{p}{2{\sigma}_4}\right)\|z^t-z^{t+1}\|^2\nonumber\\
&&-p\zeta\|x(y^{t+1}, z^t)-\bar{x}^*(z^t)\|^2.
\end{eqnarray*}
By completing the square, we further obtain
\begin{eqnarray}
%\nonumber\\
\phi^t-\phi^{t+1}
&\ge& \frac{1}{2c}\|x^{t+1}-x^t\|^2-\alpha\|A(x(y^{t+1}, z^t)-x^t)\|^2
%\nonumber\\
%&&
+\alpha\|Ax(y^{t+1}, z^t)-b\|^2\nonumber\\
&&+\left(\frac{p}{2\beta}+p-\frac{p}{\sigma_4}-\frac{p}{\zeta}-\frac{p}{2{\sigma}_4}\right)\|z^t-z^{t+1}\|^2
%\nonumber\\
%&&
-p\zeta\|x(y^{t+1}, z^t)-\bar{x}^*(z^t)\|^2\nonumber\\[3pt]
&\stackrel{ \mbox{\scriptsize(ii)}}
=& \left(\frac{1}{2c}-\frac{\alpha\smax^2(A)}{c^2(p-L_f)^2}\right)\|x^t-x^{t+1}\|^2+\alpha\|Ax(y^{t+1}, z^t)-b\|^2
\nonumber\\
&&
+\left(\frac{p}{2\beta}+p-\frac{p}{\sigma_4}-\frac{p}{\zeta}-\frac{p}{2{\sigma}_4}\right)\|z^t-z^{t+1}\|^2
%\nonumber\\
%&&
-p\zeta\|x(y^{t+1}, z^t)-\bar{x}^*(z^t)\|^2, \label{one}
\end{eqnarray}
where (i) is due to the error bound \eqref{eb1} in Lemma~\ref{error bound} and equality (ii) is due to the definition of $\sigma_1$ in Lemma~\ref{error bound}.
If we let
\begin{eqnarray}\label{alphabound}
\alpha&<&\alpha'
%\\
%&=&\frac{c^2(p-L_f)^2}{4c\smax^2(A)}\\
=\frac{c(p-L_f)^2}{4\smax^2(A)},
\end{eqnarray}
then
\begin{equation}\label{co1}
\frac{1}{2c}-\frac{\alpha\smax^2(A)}{c^2(p-L_f)^2}>\frac{1}{4c}.
\end{equation}

Moreover, according to the definition of $\sigma_4$ in Lemma \ref{error bound}, since $p\ge 3L_f$, we have
$$\sigma_4>1/2.$$
Then letting  $\zeta=6\beta$ and combining the  condition  that $\beta<1/30$, we have
\begin{equation}
p/\zeta=p/(6\beta),\quad p-p/\sigma_4-p/(2\sigma_4)\ge p/(6\beta).
\end{equation}
Therefore,
\begin{equation}\label{co2}
\left(\frac{p}{2\beta}+p-\frac{p}{\sigma_4}-\frac{p}{\zeta}-\frac{p}{2{\sigma}_4}\right)\ge (1/2-1/6-1/6)p/\beta\ge p/(6\beta).
\end{equation}
Combining \eqref{one}, \eqref{co1} and \eqref{co2}, we finish the proof of Lemma~\ref{four terms}.
\endprf

\bibliographystyle{plain}
\bibliography{reference}
\end{document}